\documentclass[11pt]{amsart}

\usepackage[headings]{fullpage}
\usepackage{amsmath}
\usepackage{amssymb,latexsym}
\usepackage{amsthm}
\usepackage{amscd}
\usepackage{bbold}

\usepackage{eucal}
\usepackage{graphicx}
\usepackage{hyperref}
\usepackage[usenames,dvipsnames]{color}
\usepackage{pdfcolmk}

\usepackage[matrix,cmtip,arrow]{xy}
\usepackage{pb-diagram,pb-xy}

\usepackage{algorithm}
\usepackage[noend]{algorithmic}
\usepackage{eqparbox}

\usepackage{xspace}

\usepackage{multirow}

\usepackage{doc}

\numberwithin{equation}{section}

\theoremstyle{plain}
\newtheorem{Theorem}[equation]{Theorem}

\newtheorem{Proposition}[equation]{Proposition}
\newtheorem{Corollary}[equation]{Corollary}

\newtheorem{theorem}[equation]{Theorem}

\newtheorem{claim}[equation]{Claim}

\theoremstyle{definition}

\theoremstyle{remark}
\newtheorem*{Remark}{Remark}
\newtheorem*{Example*}{Example}
\newtheorem*{Notation}{Notation}

\newcommand{\Hgr}                       {{\sf H}}
\newcommand{\Cgr}                       {{\sf C}}

\newcommand{\kspan}[1]{\langle{#1}\rangle}

\newcommand{\Pers}       	{\operatorname{Pers}}
\newcommand{\PI}[2]	{\left[ {#1}, {#2} \right]}
\newcommand{\Pairs}	{\operatorname{Pairs}}
\newcommand{\pair}[2]	{\left[ {#1}, {#2} \right)}

\newcommand{\Hom}	{\operatorname{Hom}}
\newcommand{\Ker}  	{\operatorname{Ker}}
\newcommand{\Img} 	 {\operatorname{Im}}

\newcommand{\td}	{\ensuremath{\circ}}		
\newcommand{\kd}	{\ensuremath{\dagger}}	

\newcommand{\ringfont}{\mathbf}
\newcommand{\kk} 	{\ringfont{k}}

\newcommand{\Zz} 	{\ringfont{Z}}

\newcommand{\ct}[1]		{{#1}^{\bot}}

\newcommand{\iso}		{\equiv}

\newcommand{\im}         	{\operatorname{im}}

\newcommand{\bdry}        	{\partial}
\newcommand{\cobdry}   	{\delta}

\newcommand{\pers}[1]	{\mathbb{#1}}

\newcommand{\Cc}		{\pers{C}}

\newcommand{\Xx}		{\pers{X}}
\newcommand{\Ss}		{\pers{S}}

\newcommand{\ssx}      	{\sigma}

\newcommand{\low}            	{\operatorname{low}}

\newcommand{\algvar}[1]    	    {{\sf\small #1}}

\newcommand{\ds}[1]             {{\sf #1}}

\newcommand{\Cpp}               {{C\hspace{-.05em}\raisebox{.4ex}{\tiny\bf ++}}\xspace}

\newcommand{\ot}		{\mathrel{\leftarrow}}

\newcommand{\phcol}{\textsf{pHcol}}
\newcommand{\phrow}{\textsf{pHrow}}
\newcommand{\pco}{\textsf{pCoh}}

\newcommand\vds[1]{}
\newcommand\dm[1]{}
\newcommand\mvj[1]{}

\title{Dualities in Persistent (Co)Homology}
\author{Vin de Silva \and Dmitriy Morozov \and Mikael Vejdemo-Johansson}
\date{\today}
\thanks{VdS has been partially supported by DARPA, through grants
  HR0011-05-1-0007 (TDA) and HR0011-07-1-0002 (SToMP), and
  holds a Digiteo Chair.
  DM has been partially supported by DARPA grant HR0011-05-1-0007 (TDA) and 
  by the DOE Office of Science, Advanced Scientific Computing Research, under
  award number KJ0402-KRD047, under contract number DE-AC02-05CH11231.
  MVJ has been partially supported by the Office of Naval Research, through grant N00014-08-1-0931.}

\begin{document}
\maketitle

\begin{abstract}
    We consider sequences of absolute and relative homology and cohomology
    groups that arise naturally for a filtered cell complex. We establish
    algebraic relationships between their persistence modules, and show that
    they contain equivalent information. We explain how one can use the existing
    algorithm for persistent homology to process any of the four modules, and
    relate it to a recently introduced persistent cohomology algorithm. We
    present experimental evidence for the practical efficiency of the latter
    algorithm.
\end{abstract}


\section{Introduction}
\label{sec:intro}

The subject of inverse problems deals, fundamentally, with the inference of shape.
From some related measurements --- such as a family of particular path integrals --- we try to deduce geometric information.
With the classical techniques in the field, with Fourier and other integral transforms,
one can deduce an impressive amount of information.
However, with non-linearity, and ill-posed, ill-conditioned situations, the classical methods need increasingly large amounts of regularization or data cleaning.
Topology offers a family of methods that allow the inference of information --- if not geometric, then at least topological --- into the field.
In particular, the recent development of persistent homology~\cite{ELZ02}, and its applications to topological data analysis~\cite{C09}, demonstrate an approach to topological invariants that becomes applicable to high-dimensional, finite and discrete measurement sets.

To take an explicit example, geological \emph{sonar} investigations employ inverse problem methods to investigate the geometric structure of the density sublevel sets in subterranean domains, relating density variations to occurrences of oil, water or mineral pockets.
The kind of information sought starts out with a qualitative judgement: is there a pocket at all; are there several or few; are they connected or not?
These first questions, before the shape can be given an explicit geometric description, are a matter of topological properties, and the study of sublevel sets of functions on domains is one of the most convincing uses of persistent homology.

The persistent homology algorithm of Edelsbrunner, Letscher, and Zomorodian~\cite{ELZ02} is now ten years old. In its natural general form~\cite{ZC05}, the input is a filtered `space' (topological space, or simplicial complex, or abstract chain complex) and the output is a collection of half-open real intervals known as a barcode or a persistence diagram.

These barcodes contain one bar for each topological feature found -- one bar for each homology class, representing a hole or a higher-dimensional void.
These bars come with a starting point, indicating the focal level at which the feature first becomes visible, and an ending point, indicating the focal level at which the feature vanishes again.
A fundamental tenet, as described in~\cite{C09} is that the length of such a bar -- the difference between when it shows up and when it vanishes -- encodes the relevance of the feature.
This emphasizes the topological features that are enveloped by a dense distribution of points, and yet have a geometrically large void in the middle. 

In many applications, all that is required is the barcode. This tells us how many homological features exist at any given level of the filtration, and how many of those survive to any given subsequent level. This information is already very rich, and has been proven to be statistically robust \cite{stability1,stability2}.
Sometimes more is required. The most common request is for geometric representatives of the features: in other words, explicit homology cycles representing each barcode interval. The original algorithm provides these cycles automatically: they are essential to the way in which the barcode intervals are calculated.

In fact, there are at least four natural persistent objects that can be derived from a filtered space. They are:
\[
\mathrm{persistent}
\left\{
\begin{array}{c}
\mathrm{absolute} \\ \mathrm{relative}
\end{array}
\right\}
\left\{
\begin{array}{c}
\mathrm{homology} \\ \mathrm{cohomology}
\end{array}
\right\}
\]
The `standard' object is persistent absolute homology, and most treatments focus on this. However, it has become increasingly clear that the other three objects are important in their own right. The transition between homology and cohomology is in some sense nothing more than the duality of vector spaces; persistent homology and cohomology have the same barcodes. However, homology cycles and cohomology cocycles are quite different, and some applications call for cocycles rather than cycles~\cite{circular}. The occasional utility of relative rather than absolute homology is probably easier to grasp intuitively; for example see~\cite{sensorcoverage} for an application in sensor networks. It is easy to `fake' the calculation of relative homology using absolute homology and a cone construction, but we point out that this trick is unnecessary.

Our goal in this paper is to provide a streamlined approach to calculating barcodes and (co)cycle representatives for all four persistent objects. We discuss this approach in terms of abstract algebra and in terms of matrix computations.

We observe that:
\begin{itemize}
\item
absolute homology and cohomology have the same barcode;
\item
relative homology and cohomology have the same barcode;
\item
the absolute barcode and the relative barcode can be deduced from each other;
\item
the cycles and bounding chains of persistent absolute homology determine, and are determined by, the cycles and bounding chains of persistent relative homology;
\item
likewise, for absolute and relative cohomology cocycles and bounding cochains.
\end{itemize}

We discuss two different dualities. There is the standard duality which interchanges homology and cohomology. We call this `pointwise' duality. More interestingly, there is a different duality which makes the following interchange:
\begin{align*}
\textrm{absolute homology} 
&\leftrightarrow 
\textrm{relative cohomology}
\\
\textrm{absolute cohomology}
&\leftrightarrow 
\textrm{relative homology}
\end{align*}
We call this `global' duality, and it appears only in the context of persistent topology. Global duality `commutes' with all possible algorithms and theorems: a method for calculating persistent absolute homology will equally well calculate persistent relative cohomology, once the input data have been turned upside-down in a particular way.

Combining all of these equalities and dualities, it emerges that a single calculation (run twice) suffices to calculate all four persistent objects. Actually, we describe two different algorithms for that calculation: \phcol\ (the `column algorithm') and \phrow\ (the `row algorithm'). Here \phcol\ is essentially the classic algorithm of~\cite{ELZ02,ZC05}; \phrow\ organises the calculation quite differently. The preferred choice depends, in any given situation, on whether it is easier to look up rows or columns of the boundary matrix of the filtered space --- the specific representation of the space usually biases this choice.

We are rewarded by an unexpected payoff. If we require only the absolute barcode, it turns out that the best choice is an optimised version of \phrow\, called \pco\ (the `cohomology algorithm'). We give experimental evidence to this effect. Standard practice has been to use \phcol. We therefore call on persistent topology library-writers to implement \pco, and on persistent topology library-users to use it.

\subsection{Outline of paper}
\label{sec:outline-paper}

Section~\ref{sec:algebra} is devoted to the algebra underlying this work. In \ref{subsec:coeffs}--\ref{subsec:stand-barc-isom} we conduct the discussion at a high level (homology functors are assumed given, black-box style), and in \ref{subsec:chains}--\ref{subsec:cohomology} we go into the necessary chain-level details. In \ref{subsec:global} we give a brief abstract description of the two dualities.

Section~\ref{sec:matrix-presentation} is about matrix algorithms. In \ref{sec:interpreting-d}--\ref{sec:cohogen} we interpret the preceding algebra in terms of matrix decomposition (following~\cite{vineyards}, again black-box style). In~\ref{sec:elz} we present the two algorithms, \phcol\ and \phrow, and explain why they give the same output.

In Section~\ref{sec:optimizations} we relate the ideas in this paper to an earlier cohomology algorithm \pco\ published in~\cite{circular}. We indicate why we expect \pco\ to be faster that \phcol\ and \phrow\ for computing barcodes of filtered simplicial complexes, and we verify this by experiment.

\section{Algebra}
\label{sec:algebra}

We will assume that the reader is familiar with homology theory. Our preference is to use cellular homology, because it is a little more general than simplicial homology.

\subsection{Coefficients}
\label{subsec:coeffs}
Individual (co)homology groups are defined with coefficients in a field~$\kk$, which remains fixed throughout this paper. Persistent (co)homology then has the structure of a graded module over the polynomial ring $\kk[t]$. Many things go wrong when we replace the field~$\kk$ with a ring, in particular the ring of integers $\Zz$. See~\cite{ZC05}.

\subsection{Filtered complexes}
We are interested in the persistent topology of filtered topological spaces. The simplest example is a filtered cell complex, which is a sequence $\Xx$ of cell complexes
\begin{equation}
\Xx: \quad
X_1 \subset X_2 \subset \dots \subset X_n = X_\infty
\label{eqn:filtcomplex}
\end{equation}
where $X_1$ is a vertex~$\sigma_1$, and thereafter each complex is obtained from the previous one by adding a single cell: $X_{i} = X_{i-1} \cup \sigma_i$.
Here the index set is $\{1, 2, \dots, n\}$. Usually we attach real values $a_i$ to the indices, which must satisfy $a_1 \leq a_2 \leq \dots \leq a_n$.

\begin{Example*}
Our running example $\Ss$ will be a cellular filtration of the 2-sphere:
\begin{center}
\includegraphics[scale=0.7,page=2]{figures/filtration}
\includegraphics[scale=0.7,page=3]{figures/filtration}
\includegraphics[scale=0.7,page=4]{figures/filtration}
\includegraphics[scale=0.7,page=5]{figures/filtration}
\includegraphics[scale=0.7,page=6]{figures/filtration}
\includegraphics[scale=0.7,page=7]{figures/filtration}
\end{center}
There are six cells, $\sigma_1, \dots, \sigma_6$ which appear at times $a_i = i$, for $i = 1, \dots, 6$.
\end{Example*}

\subsection{Persistent homology}

If we apply a homology functor $\Hgr(-)$ to a filtered complex~$\Xx$ we obtain a diagram:
\begin{equation}
\Hgr(\Xx):\quad
\Hgr(X_1) \to
\Hgr(X_2) \to
\dots
\to
\Hgr(X_n)
\label{eqn:Hfiltcomplex}
\end{equation}
Typically $\Hgr(-)$ denotes the $k$-dimensional homology $\Hgr_k(-; \kk)$ or the total homology $\Hgr_*(-; \kk)$. Then \eqref{eqn:Hfiltcomplex} is a diagram of finite-dimensional vector spaces and linear maps, also known as a \textbf{persistence module}.

A persistence module decomposes as a direct sum of \textbf{interval modules}~\cite{ZC05}. These are labelled by ordered pairs of integers $\PI{p}{q}$, where $1 \leq p \leq q \leq n$. The pair $\PI{p}{q}$ indicates a feature which persists over the index set $\{p, \dots, q\}$. We frequently interpret $\PI{p}{q}$ as the half-open real interval $[a_p, a_{q+1})$, with the convention that $a_{n+1} = \infty$.

The \textbf{persistence diagram} or \textbf{barcode} is the multiset of ordered pairs $\PI{p}{q}$ in the decomposition, or alternatively the multiset of half-open intervals $[a_p, a_{q+1})$. Thus we write:
\begin{align*}
\Pers(\Hgr(\Xx))
&= \left\{ \PI{p_1}{q_1}, \dots, \PI{p_m}{q_m} \right\}
\\
&= \left\{
	\left[ a_{p_1}, a_{q_1+1} \right), \dots,
	\left[ a_{p_m}, a_{q_m+1} \right)
      \right\}
\end{align*}
It is customary in applications to discard from the persistence diagram those intervals $[a_p, a_{q+1})$ for which $a_p = a_{q+1}$.

\begin{Example*}
In our running example, the intermediate spaces $S_1, S_3, S_5$ are all contractible, whereas $S_2, S_4, S_6$ are homeomorphic to the 0-sphere, 1-sphere, and 2-sphere, respectively. There are four intervals in the persistence diagram of $\Hgr_*(\Ss)$:
\[
\Pers(\Hgr_*(\Ss))
= \left\{ \PI{1}{6}_0, \PI{2}{2}_0, \PI{4}{4}_1, \PI{6}{6}_2 \right\}
= \left\{
	[1,\infty)_0, [2,3)_0, [4,5)_1, [6,\infty)_2
      \right\}
\]
The subscript~$k$ in $\PI{p}{q}_k$ or $[a_p, a_{q+1})_k$ indicates that the feature occurs in $k$-dimensional homology.
\end{Example*}

\subsection{The four standard persistence modules}
\label{sec:four-stand-pers}

The standard persistent homology module $\Hgr_*(\Xx)$ tells us how the absolute homology groups $\Hgr_*(X_i)$ relate to each other as $i$~varies. We can play the same game with the absolute cohomology groups $\Hgr^*(X_i)$, the relative homology groups $\Hgr_*(X_n, X_{i})$, and the relative cohomology groups $\Hgr^*(X_n, X_i)$. Here are the four sequences, lined up for comparison.
\begin{alignat*}{5}
\Hgr_*(\Xx) :\quad
&
&& \mathrel{\phantom{\to}}
\Hgr_*(X_1) && \to {\;\dots\;} && \to \Hgr_*(X_{n-1}) && \to \Hgr_*(X_{n})
\\
\Hgr^*(\Xx) :\quad
&
&& \mathrel{\phantom{\ot}}
\Hgr^*(X_1) && \ot {\;\dots\;} && \ot \Hgr^*(X_{n-1}) && \ot \Hgr^*(X_{n})
\\
\Hgr_*(X_\infty, \Xx) :\quad
& \Hgr_*(X_n) && \to \Hgr_*(X_n, X_1) && \to {\;\dots\;} && \to \Hgr_*(X_n, X_{n-1})
\\
\Hgr^*(X_\infty, \Xx) :\quad
& \Hgr^*(X_n) && \ot \Hgr^*(X_n, X_1) && \ot {\;\dots\;} && \ot \Hgr^*(X_n, X_{n-1})
\end{alignat*}

The persistence diagram for absolute cohomology is a multiset of integer ordered pairs $\PI{p}{q}$ with $1 \leq p \leq q \leq n$. For relative homology and cohomology, the persistence diagrams are multisets of pairs $\PI{p}{q}$ with $0 \leq p \leq q \leq n-1$. In all cases, we interpret $\PI{p}{q}$ as the half-open interval $[a_p, a_{q+1})$, with the convention that $a_0 = -\infty$ and $a_{n+1} = \infty$.

\begin{Example*}
In our running example, we compute
\[
\Pers(\Hgr_*(S_6, \Ss))
= \left\{ \PI{0}{0}_0, \PI{2}{2}_1, \PI{4}{4}_2, \PI{0}{5}_2 \right\}
= \left\{
	[-\infty, 1)_0, [2,3)_1, [4,5)_2, [-\infty,6)_2
      \right\}.
\]
For instance, at index~2 we note that there is a nontrivial element of $\Hgr_1(S_6, S_2)$ represented by any arc connecting the two points of $S_2$. To be specific, the homology class is $[\sigma_3] = [\sigma_4]$. This class vanishes in $\Hgr_1(S_6,S_3)$, and so it generates the interval $[2,3)$.
\end{Example*}

The reader may detect a relationship between the barcodes for absolute and relative homology. We formalize this in the next section.

\subsection{Barcode isomorphisms}
\label{subsec:stand-barc-isom}

\begin{Proposition}
\label{prop:UCTk}
For all~$k$,
\begin{eqnarray*}
\Pers(\Hgr_k(\Xx)) &=& \Pers(\Hgr^k(\Xx)),
\\
\Pers(\Hgr_k(X_\infty, \Xx)) &=& \Pers(\Hgr^k(X_\infty, \Xx)).
\end{eqnarray*}
In other words, homology and cohomology have identical barcodes.
\end{Proposition}

\begin{proof}
The universal coefficients theorem \cite[Thm 3.2]{hatcher} asserts that there is a natural isomorphism
\[
\Hgr^k(X; \kk) \iso \Hom(\Hgr_k(X; \kk), \kk).
\]
In other words, cohomology and homology are dual as vector spaces, and hence have the same dimension. `Natural' implies that the induced maps
\[
\Hgr_k(X_i; \kk) \to \Hgr_k(X_j; \kk)
\quad \mbox{and} \quad
\Hgr^k(X_i; \kk) \ot \Hgr^k(X_j; \kk)
\]
are adjoint, and hence have the same rank. Because of the way the barcode is uniquely determined by dimensions and ranks, it follows that the absolute homology and cohomology barcodes are the same. This argument applies equally well to the relative barcodes.
\end{proof}

\begin{Notation}
We partition each persistence diagram into two parts,
\[
\Pers = \Pers_0 \cup \Pers_\infty,
\]
where $\Pers_0$ comprises the finite intervals $[a,b)$, and $\Pers_\infty$ the infinite intervals $[a,\infty)$ or $[-\infty, b)$.
\end{Notation}

\begin{Proposition}
\label{prop:UCTkt}
For all~$k$,
\begin{eqnarray*}
\Pers_0(\Hgr_k(\Xx)) &=&
 \Pers_0(\Hgr_{k+1}(X_\infty, \Xx)),
\\
\Pers_\infty(\Hgr_k(\Xx)) &=&
 \Pers_\infty(\Hgr_{k}(X_\infty, \Xx)),
\end{eqnarray*}
where the second `equality' is interpreted as a bijection with $[a,\infty) \leftrightarrow [-\infty, a)$.
Thus, persistent homology and relative homology barcodes carry the same information, with a dimension shift for the finite intervals.
\end{Proposition}

The proof appears in Section~\ref{subsec:chains}.


\begin{Remark}
Thus, provided we take the dimension shifts into account, all four barcodes carry exactly the same information. If we are only interested in barcodes, we can perform calculations in any one of the four basic sequences, whichever is the most convenient.
\end{Remark}


\vds{This next passage is optional, although it does shed light on why extended persistence relativizes the space in reverse direction.}

Since the last term of $\Hgr_k(\Xx)$ is the same as the first term of $\Hgr_k(X_\infty, \Xx)$, namely $\Hgr_k(X_n)$, the two sequences can be concatenated into a single sequence, which we denote $\Hgr_k(\Xx) \to \Hgr_k(X_\infty, \Xx)$. The index set for this sequence is $\{ 1, 2, \dots, n = \bar{0}, \bar{1}, \bar{2}, \dots, \overline{n-1}\}$, where we use barred numerals to indicate that we are in the relative homology part of the sequence. The persistence diagram for this complex will have intervals of three possible types:
\begin{itemize}
\item
$(p,q)$ where $1 \leq p \leq q < n$, written as $[p,q+1)$ or $[a_p, a_{q+1})$ in interval form.
\item
$(\bar{p}, \bar{q})$ where $0 < p \leq q \leq n-1$, written as $[\bar{p}, \overline{q+1})$ or $[\bar{a}_p, \bar{a}_{q+1})$.
\item
$(p,\bar{q})$ where $1 \leq p \leq n$, $0 \leq q \leq n-1$, written as $[p, \overline{q+1})$ or $[a_p, \bar{a}_{q+1})$.
\end{itemize}

\begin{Proposition}
\label{prop:abs-rel}
The barcode $\Pers\left(\Hgr_k(\Xx) \to \Hgr_k(X_\infty, \Xx)\right)$ comprises the following collection of intervals:
\begin{itemize}
\item
An interval $[a,b)$ for every interval $[a,b)$ in $\Pers_0(\Hgr_k(\Xx))$.
\item
An interval $[\bar{a},\bar{b})$ for every interval $[a,b)$ in $\Pers_0(\Hgr_{k-1}(\Xx))$.
\item
An interval $[a, \bar{a})$ for every interval $[a, \infty)$ in $\Pers_\infty(\Hgr_k(\Xx))$.
\end{itemize}
\end{Proposition}

\begin{proof}
Note that the first two classes of interval in $\Pers\left( \Hgr_k(\Xx) \to \Hgr_k(X_\infty, \Xx) \right)$ are those which do not meet the middle term $\Hgr_k(X_n)$, and thus correspond exactly to finite intervals in $\Pers(\Hgr_k(\Xx))$ and $\Pers(\Hgr_k(X_\infty, \Xx))$. This explains the first two cases, once we make the translation $\Pers_0(\Hgr_k(X_\infty, \Xx)) = \Pers_0(\Hgr_{k-1}(\Xx))$.

It remains to show is that the intervals of type $[a, \bar{b})$ are always of the form $[a, \bar{a})$.\footnote{%
Thus the paired intervals $[a, \infty)$ and $[-\infty, a)$ in $\Pers_\infty(\Hgr_k(\Xx))$ and $\Pers_\infty(\Hgr_k(X_\infty, \Xx))$ are really the restrictions of a single interval $[a, \bar{a})$ in the concatenated sequence.%
}
To do this, we need to compare the right filtration of the sequence $\Hgr_k(\Xx)$ with the left filtration of the sequence $\Hgr_k(X_\infty, \Xx)$.
\vds{Need to cite zigzag paper.}
The first filtration is the nested sequence of subspaces
\[
\Img(\Hgr_k(X_i) \to \Hgr_k(X_n)),
\qquad
i = 1, 2, \dots, n-1,
\]
of $\Hgr_k(X_n)$, and the second filtration is the nested sequence of subspaces
\[
\Ker(\Hgr_k(X_n) \to \Hgr_k(X_n, X_i))
\qquad
i = 1, 2, \dots, n-1,
\]
of $\Hgr_k(X_n)$.
But the image and kernel subspaces are equal for each~$i$, by the homology long exact sequence for the pair $(X_n, X_i)$. Thus the filtrations are the same.
\end{proof}

\begin{Remark}
The sequence $\Hgr_k(\Xx) \to \Hgr_k(X_\infty, \Xx)$ is not the same as the
extended persistence~\cite{cohen2009extending}. The latter, defined for the
sublevel sets of a real-valued function, requires the reversal of the cells in
the relative half of the sequence --- it translates into the use of the
superlevel sets of the function.
The meaning of extended persistence for a general filtered space is a
lot less straight-forward. The most significant difference between the two
sequences (besides the definition) are the extended pairs, the intervals
corresponding to $[a, \infty)$ in $\Pers_\infty(\Hgr_k(\Xx))$. In
Proposition~\ref{prop:abs-rel} they become the trivial intervals $[a,
\bar{a})$; on the other hand, they are the main reason extended persistence was
introduced: these pairs carry the new information. Another notable difference
is that the dualities in this paper apply to general filtered spaces;
Poincar\'e and Lefschetz dualities involved in the analysis of the extended
persistence require the domains to be manifolds.
\end{Remark}

\vds{end optional section}

\subsection{Persistent chain complexes}
\label{subsec:chains}

We now give a more explicit description of the standard persistence modules, in terms of chain complexes. Among other things, this will lead to a clean proof of Proposition~\ref{prop:UCTkt}.
Given a filtered cell complex~$\Xx = \sigma_1 \cup \dots \cup \sigma_n$, define a persistence module
\[
\Cc: \quad
C_1 \to C_2 \to \dots \to C_n
\]
where $C_i = \kspan{\sigma_1, \dots, \sigma_i}$, the vector space over~$\kk$ with basis elements labelled $\sigma_1, \dots, \sigma_i$. We also have a boundary map: the boundary of any $\sigma_j$ is a linear combination of cells which appear previously:
\[
\bdry \sigma_j = \sum_{i < j} {D}_{ij} \sigma_i
\]
for some collection of coefficients $D_{ij}$.
Geometrically, the cells $\sigma_i$ for which $D_{ij} \ne 0$ will have dimension one less than the dimension of~$\sigma_j$.

The boundary map satisfies $\bdry^2 = 0$, and it restricts to boundary maps $\bdry_i: C_i \to C_i$ for each~$i$. Then $\Cgr_*(X_i) = (C_i, \bdry_i)$ is the chain complex%
\footnote{%
For simplicity we generally suppress the homological grading within each $C_i$, which comes from the geometric dimensions of the cells associated to the generators. We will refer to it only when necessary.%
}
for the absolute homology of~$X_i$, and $\Cgr_*(\Xx) = (\Cc, \bdry)$ is the persistent version for $\Xx$. Accordingly, we define the persistent absolute homology of~$\Xx$ to be
\[
\Hgr_*(\Xx) = \Hgr(\Cc, \bdry) = \Ker(\bdry) / \Img(\bdry).
\]
Shown explicitly as a persistence module, this is:
\[
\Hgr(\Cc, \bdry):\quad
\frac{\Ker(\bdry_1)}{\Img(\bdry_1)}
\to
\frac{\Ker(\bdry_2)}{\Img(\bdry_2)}
\to
\dots
\to
\frac{\Ker(\bdry_n)}{\Img(\bdry_n)}.
\]

For the absolute cohomology persistence module $\Hgr^*(\Xx)$, we define
\[
\Cc^*: \quad
C_1^* \ot C_2^* \ot \dots \ot C_n^*
\]
where $C_i^* = \Hom(C_i, \kk) = \kspan{\sigma_1^*, \sigma_2^*, \dots, \sigma_i^*}$, with $\{\sigma_i^*\}$ being the dual basis to $\{\sigma_i\}$. The coboundary $\cobdry = \bdry^*$ is defined to be the adjoint to $\bdry$. Then $\Cgr^*(\Xx) = (\Cc^*, \cobdry)$ and
\[
\Hgr^*(\Xx) = \Hgr(\Cc^*, \cobdry) = \Ker(\cobdry) / \Img(\cobdry).
\]
Again, this is a persistence module (with arrows to the left).

\begin{Example*}
In our running example, the boundary map is given as follows:
\begin{align*}
&\bdry \sigma_1 = \bdry \sigma_2 = 0,
\\
&\bdry\sigma_3 = \bdry\sigma_4 = \sigma_1 - \sigma_2,
\\
&\bdry\sigma_5 = \bdry\sigma_6 = \sigma_3 - \sigma_4.
\end{align*}
This information is recorded in matrix~$D$ of Figure~\ref{fig:matrix-example}.
The coboundary map is given as follows:
\begin{align*}
&\cobdry \sigma^*_1 = - \cobdry\sigma_2 = \sigma^*_3 + \sigma^*_4,\;
\\
&\cobdry \sigma^*_3 = -\cobdry\sigma^*_4 = \sigma^*_5 + \sigma^*_6,\;
\\
&\cobdry \sigma^*_5 = \phantom{-}\cobdry\sigma^*_6 = 0.
\end{align*}
This information is recorded in matrix~$\ct{D}$ of Figure~\ref{fig:matrix-example}.
\end{Example*}

The relative homology and cohomology persistence modules are defined as the homology of the persistence modules
\begin{alignat*}{5}
(C_n/\Cc)^{\phantom{*}}
&:\quad
C_n
&& \to
(C_n/C_1)
&& \to
(C_n/C_2)
&& \to
\;\dots\;
&& \to
(C_n/C_{n-1})
\\
(C_n/\Cc)^* &:\quad
C_n^*
&& \ot
(C_n/C_1)^*
&& \ot
(C_n/C_2)^*
&& \ot
\;\dots\;
&& \ot
(C_n/C_{n-1})^*
\end{alignat*}
with boundary and coboundary maps induced from $\bdry, \cobdry$ in the natural way. Thus
\begin{align*}
\Cgr_*(X_\infty, \Xx) &= (C_n/\Cc, \bdry),
&
\Cgr^*(X_\infty, \Xx) &= ((C_n/\Cc)^*, \cobdry).
\\
\Hgr_*(X_\infty, \Xx) &= \Hgr(C_n/\Cc, \bdry),
&
\Hgr^*(X_\infty, \Xx) &= \Hgr((C_n/\Cc)^*, \cobdry).
\end{align*}

\begin{Remark}
We note that the maps $\to$ of $\Cc$ and the maps $\ot$ of $(C_n / \Cc)^*$ are injective, whereas the maps $\ot$ of $\Cc^*$ and the maps $\to$ of $(C_n/\Cc)$ are surjective. In other words, absolute homology and relative cohomology are structurally akin to each other; and qualitatively different from absolute cohomology and relative homology. This is a symptom of the global duality mentioned in the introduction.
\end{Remark}

The main theorem of~\cite{ELZ02,ZC05} can be restated in the following way.

\begin{Theorem}
\label{thm:persistence-main}
Given $\Cc, \bdry$ as above, there exists a partition
\[
\{ 1, 2, \dots, n \} = F \sqcup G \sqcup H
\]
with a bijective pairing $G \leftrightarrow H$, written as follows:
\[
\mbox{$g$ is paired with $h$}
\; \Leftrightarrow \;
\pair{g}{h} \in \Pairs = \Pairs(\Cc, \bdry).
\]
Moreover, there is a new basis $\hat\sigma_1, \hat\sigma_2, \dots, \hat\sigma_n$ of $C_n$ such that:
\begin{enumerate}
\item
$C_i = \kspan{\hat\sigma_1, \dots, \hat\sigma_i}$ for all~$i$.
\item
$\bdry \hat\sigma_f = 0$ for all $f \in F$.
\item
$\bdry \hat\sigma_h = \hat\sigma_g$, and hence $\bdry \hat\sigma_g = 0$, for all $\pair{g}{h} \in \Pairs$.
\end{enumerate}
It follows that the persistence diagram $\Pers(\Hgr(\Cc,\bdry))$ consists of the intervals $[a_f,\infty)$ for $f \in F$ together with the intervals $[a_g,a_h)$ for $[g,h) \in \Pairs$.
\qed
\end{Theorem}

We note that item~(1) is equivalent to the assertion that the leading term of each $\hat\sigma_i$ is $\sigma_i$ (up to a nonzero scalar).

In the language of~\cite{ELZ02}, the index set~$F$ identifies the positive simplices which remain unpaired, the index set~$G$ identifies the positive simplices which do get paired, and the index set~$H$ identifies the negative simplices. The vectors $\hat\sigma_f$ and~$\hat\sigma_g$ are the cycles with leading terms $\sigma_f$ and~$\sigma_g$, and the vector~$\hat\sigma_h$ is the chain with leading term~$\sigma_h$ which `kills' the homology class of its paired~$\hat\sigma_g$ by means of the equation $\bdry\hat\sigma_h = \hat\sigma_g$.

%

\begin{Example*}
In our running example, $F = \{1,6\}$, $G = \{2, 4\}$, $H = \{3,5\}$ and $\Pairs = \{ \pair{2}{3}, \pair{4}{5} \}$.
The new basis is
\begin{align*}
\hat\sigma_1 &= \sigma_1,
&
\hat\sigma_3 &= \sigma_3,
&
\hat\sigma_5 &= \sigma_5,
\\
\hat\sigma_2 &= -\sigma_2 + \sigma_1,
&\hat\sigma_4 &= -\sigma_4 + \sigma_3,
&\hat\sigma_6 &= \sigma_6 - \sigma_5.
\end{align*}
The reader can easily verify that $\bdry\hat\sigma_1 = \bdry\hat\sigma_2 = \bdry\hat\sigma_4 = \bdry\hat\sigma_6 = 0$, that $\bdry\hat\sigma_3 = \hat\sigma_2$, and $\bdry\hat\sigma_5 = \hat\sigma_4$.
\end{Example*}

\begin{proof}[Proof of Proposition~\ref{prop:UCTkt}]
The decomposition $\{1, 2, \dots, n\} = F \sqcup G \sqcup H$ and the new basis $\hat\sigma_1, \hat\sigma_2, \dots, \hat\sigma_n$ allow us to express $\Cc$ as a direct sum of persistent chain complexes:
\[
\Cc =
 \bigoplus_{f \in F} \Cc_f
 \quad\oplus\quad
 \bigoplus_{\pair{g}{h} \in \Pairs} \Cc_{{g},{h}}
\]
where $\Cc_f = \kspan{\hat\sigma_f}$ and $\Cc_{{g},{h}} = \kspan{\hat\sigma_g, \hat\sigma_h}$.
Moreover, the boundary map $\bdry$ respects this decomposition, mapping each summand into itself.
We can therefore calculate $\Hgr((C_n/\Cc), \bdry)$ on each summand separately.

For summands of type $\Cc_f$, the persistence modules are constant over two phases, with index ranges $\{0, \dots, f-1\}$ and $\{f, \dots, n-1\}$. We can condense this information by representing them as two-term persistence modules (one term for each index range):
%
%
\begin{alignat*}{2}
((C_f)_n/ \Cc_f)
&:\quad
\kspan{\hat\sigma_f}
&&\to 0
\\
\Ker(\bdry)
&:\quad
\kspan{\hat\sigma_f}
&&\to 0
\\
\Img(\bdry)
&:\quad
0
&&\to 0
\\
\Hgr = \Ker(\bdry)/\Img(\bdry)
&:\quad
\kspan{\hat\sigma_f}
&&\to 0
\end{alignat*}
It follows that $\Hgr((C_f)_n/ \Cc_f)$ contributes the interval $[-\infty, a_f)$. This is generated by $[\hat\sigma_f]$ and hence has the same homological dimension as $[a_f,\infty)$ in $\Pers(\Hgr(\Cc))$.

For summands of type $\Cc_{{g},{h}}$ the persistence modules are constant over three phases, with index ranges $\{0, \dots, g-1\}$,  $\{g, \dots, h-1\}$ and $\{h, \dots, n-1\}$:
\begin{alignat*}{3}
((C_{{g},{h}})_n/ \Cc_{{g},{h}})
&:\quad
\kspan{\hat\sigma_g, \hat\sigma_h}
&&\to
\kspan{\hat\sigma_h}
&&\to 0
\\
\Ker(\bdry)
&:\quad
\kspan{\hat\sigma_g}
&&\to \kspan{\hat\sigma_h}
&&\to 0
\\
\Img(\bdry)
&:\quad
\kspan{\hat\sigma_g}
&&\to 0
&&\to 0
\\
\Hgr = \Ker(\bdry)/\Img(\bdry)
&:\quad
0
&&\to
\kspan{\hat\sigma_h}
&&\to 0
\end{alignat*}
It follows that $\Hgr((C_{g,h})_n/ \Cc_{g,h})$ contributes a single interval, $[a_g, a_h)$. This is generated by $[\hat\sigma_h]$, and hence has dimension one greater than $[a_g, a_h)$ in $\Pers(\Hgr(\Cc))$, that being generated by $[\hat\sigma_g]$.
\end{proof}

 The following table summarises the relationship between the three
  types of generator and the persistence intervals they generate.
  \medskip
\begin{center}
{
\begin{tabular}{|c||c||c|c|}
\hline
generator
  & $\hat\sigma_f $
  & $\hat\sigma_g $
  & $\hat\sigma_h $
\\
\hline
absolute homology & $[a_f, \infty)$ & $[a_g, a_h)$ &
\\
\hline
relative homology & $[-\infty, a_f)$ & & $[a_g, a_h)_{_+}$
\\
\hline
\end{tabular}
}
\end{center}
\medskip
The homological dimension of each interval is equal to the homological dimension of its generator. So, for any pair $[g,h) \in \Pairs$, the dimension of $[a_g,a_h)$ in $\Hgr_*(X_\infty, \Xx)$ is one greater than the dimension of $[a_g,a_h)$ in $\Hgr_*(\Xx)$. We indicate this in the table with a $_+$~subscript.

\subsection{Cohomology}
\label{subsec:cohomology}

Persistent relative cohomology is structurally similar to persistent absolute homology. To make this apparent, let us introduce new notation, writing 
\begin{alignat*}{4}
\ct{\Cc}
&:\quad
\ct{C}_{1}
&& \to
\ct{C}_{2}
&& \to
\;\dots\;
&& \to
\ct{C}_{n}
\end{alignat*}
for the reverse of the sequence
\begin{alignat*}{4}
(C_n/\Cc)^* &:\quad
C_n^*
&&\ot
(C_n/C_1)^*
&& \ot
\;\dots\;
&& \ot
(C_n/C_{n-1})^*,
\end{alignat*}
so $\ct{C}_i = (C_n / C_{n-i})^*$.

Recall that $\sigma_1^*, \dots, \sigma_n^*$ denotes the basis of $C_n^*$ dual to the basis $\sigma_1, \dots, \sigma_n$ of $C_n$. If we write $\tau_i = \sigma_{n+1-i}^*$, then
\[
\ct{C}_{i} = (C_n / C_{n-i})
= \kspan{\sigma_n^*, \sigma_{n-1}^*, \dots, \sigma_{n+1-i}^*}
= \kspan{\tau_1, \tau_2, \dots, \tau_i}
\]
by elementary linear algebra. Moreover, if
$
\bdry \sigma_j = \sum_{i < j} D_{ij} \sigma_i
$
then
$
\cobdry \tau_j = \sum_{i < j} \ct{D}_{ij} \tau_i.
$
where $\ct{D}_{ij} = D_{(n-j),(n-i)}$.

On account of the formal similarity between $\ct{\Cc}$ and $\Cc$ we conclude: 

\begin{Corollary}
Suppose we have an algorithm which takes as input the sequence of cells $\sigma_1, \dots, \sigma_n$ and their boundaries $\bdry\sigma_1, \dots, \bdry\sigma_n$ and produces as output
the persistent absolute homology of $\Xx = \sigma_1 \cup \dots \cup \sigma_n$.

Then the same algorithm applied to the sequence of formal cells $\tau_1, \dots, \tau_n$ and coboundaries $\cobdry\tau_1, \dots, \cobdry\tau_n$ computes the persistent relative cohomology of~$\Xx$.
\qed
\end{Corollary}

Persistent absolute cohomology can be thought of as `relative persistent relative cohomology'.
More precisely, by elementary linear algebra:
\begin{Proposition}
The persistence module $(\ct{C}_n / \ct{\Cc})$ is the reverse of $\Cc^*$, and the respective coboundary maps agree.
\qed
\end{Proposition}

\begin{Corollary}
Suppose we have an algorithm which takes as input the sequence of cells $\sigma_1, \dots, \sigma_n$ and their boundaries $\bdry\sigma_1, \dots, \bdry\sigma_n$ and produces as output the persistent relative homology of $\Xx = \sigma_1 \cup \dots \cup \sigma_n$.

Then the same algorithm applied to the sequence of formal cells $\tau_1, \dots, \tau_n$ and coboundaries $\cobdry\tau_1, \dots, \cobdry\tau_n$ computes the persistent absolute cohomology of~$\Xx$.
\qed
\end{Corollary}

\begin{Remark}
We must transcribe the indices correctly for these two corollaries. If such an algorithm applied to the $\tau_i, \cobdry\tau_i$ produces a persistence interval $\PI{p}{q-1}$ for $\ct{\Cc}$ then this is equivalent to $\PI{n+1-q}{n-p}$ for $(C_n/ \Cc)^*$, and hence to the half-open real interval $[a_{n+1-q}, a_{n+1-p})$.
\end{Remark}

Suppose now we apply Theorem~\ref{thm:persistence-main} to $\ct{\Cc}$ to obtain a partition $\{1, 2, \dots, n\} = R \sqcup S \sqcup T$ and new generators $\hat\tau_i$, with $\cobdry\hat\tau_r = 0$, $\cobdry\hat\tau_s = 0$, and $\cobdry\hat\tau_t = \hat\tau_s$ for pairs $[s,t) \in \Pairs(\ct{\Cc},\cobdry)$.
Then we obtain the following table:
\medskip
\begin{center}
{\small
\begin{tabular}{|c||c||c|c|}
\hline
generator
  & $\hat\tau_{r} $
  & $\hat\tau_{s} $
  & $\hat\tau_{t} $
\\
\hline
relative cohomology
  & $[-\infty, a_{n+1-r})$
  & $[a_{n+1-t}, a_{n+1-s})_{_+}$
  &
\\
\hline
absolute cohomology
  & $[a_{n+1-r}, \infty)$ 
  &
  & $[a_{n+1-t}, a_{n+1-s})$
\\
\hline
\end{tabular}
}
\end{center}
\medskip
By considering Proposition~\ref{prop:UCTk}, we deduce that
\[
R = n + 1 - F, \quad
S = n + 1 - H, \quad
T = n + 1 - G
\]
and moreover $(s,t) \in \Pairs(\ct{\Cc},\cobdry)$ if and only if $(n+1-t, n+1-s) \in \Pairs(\Cc,\bdry)$. Actually, this can also be inferred from the proof of the following proposition.

\begin{Proposition}
Let $\hat\sigma_1^*, \dots,  \hat\sigma_n^*$ denote the dual basis to $\hat\sigma_1, \dots, \hat\sigma_n$, and write $\hat\tau_i = \hat\sigma_{n+1-i}^*$. Then the $\hat\tau_i$ are the generators described above (up to nonzero scalar multiples).
\end{Proposition}

\begin{proof}
By duality, we have $\cobdry \hat\sigma^*_f = 0$ for all $f \in F$, and $\cobdry \hat\sigma^*_g = \hat\sigma^*_h$ for all $\pair{g}{h} \in \Pairs(\Cc)$. Moreover, $\sigma^*_i$ is the \emph{trailing} term of $\hat\sigma_i^*$. The proposition now follows, by bookkeeping.
\end{proof}

\subsection{A remark for the algebraically-minded}
\label{subsec:global}

According to \cite{ZC05}, a persistence module can be regarded as a graded module over the ring $\kk[t]$. In particular, $\Cc$ can  be regarded as a free module over $\kk[t]$ with $n$ generators, $\sigma_1, \dots, \sigma_n$, where $\sigma_i$ has degree~$i$. The boundary map $\bdry: \Cc \to \Cc$ is then a homomorphism of graded modules.

The `dual' of such a module can be taken with respect to the ground field $\kk$ or the polynomial ring $\kk[t]$. Thus we can define the global dual
\begin{align*}
\Cc^\td &= \Hom_{\kk[t]}(\Cc, \kk[t]),
\\
\mbox{where}\quad
C^\td_n &=
\makebox[23em][l]{graded-module homomorphisms $\Cc \to \kk[t]$ of degree~$n$;}
\end{align*}
and the pointwise dual
\begin{align*}
\Cc^\kd &= \Hom_\kk(\Cc, \kk),
\\
\mbox{where}\quad
C^\kd_n &=
\makebox[23em][l]{vector-space homomorphisms $C_{-n} \to \kk$.}
\end{align*}
These can be regarded as graded $\kk[t]$-modules in a natural way. Moreover, the operations $-^\td$ and $-^\kd$ are contravariant functors, so in particular the boundary map on $\Cc$ induces boundary maps on the new modules.

The interested reader can verify that
\begin{align*}
\Hgr(\Cc, \bdry) &= \mbox{persistent absolute homology of $\Xx$}
\\
\Hgr(\Cc^\kd, \bdry^\kd) &= \mbox{persistent absolute cohomology of $\Xx$}
\\
\phantom{^\kd}
\Hgr(\Cc^\td, \bdry^\td) &= \mbox{persistent relative cohomology of $\Xx$}
\\
\Hgr(\Cc^{\td\kd}, \bdry^{\td\kd}) &= \mbox{persistent relative homology of $\Xx$}
\end{align*}
up to calibrating the indices.

\section{Matrix Algorithms}
\label{sec:matrix-presentation}

\subsection{The boundary matrix}
\label{sec:interpreting-d}

We can represent a filtered cell complex (at least, its homological information) by a strictly upper-triangular matrix~$D$, whose entries $D[i,j]$ are the coefficients $D_{ij}$ of the boundary map~$\bdry$ defined in Section~\ref{subsec:chains}. Thus the $j$-th column $D[j] = D[..,j]$ represents $\bdry \sigma_j$.
With the cells listed in the filtration order, the matrix~$D$ also encodes the filtration of the complex. Indeed, the top-left square submatrix $D[1..i,1..i]$ is the boundary matrix for $X_i = \sigma_1 \cup \dots \cup \sigma_i$, or, equivalently, the chain complex $C_i$. Thus $D$ is a representation of the chain complex for persistent absolute homology, $\Cgr_*(\Xx) = (\Cc, \bdry)$.

If we flip the matrix~$D$ across its minor diagonal, we get the \textbf{anti-transpose} $\ct{D}$, formally defined by $\ct{D}[i,j] = D[n+1-j,n+1-i]$. Following the discussion in Section~\ref{subsec:cohomology}, we see that $\ct{D}$ represents the cochain complex for persistent relative cohomology, $\Cgr^*(X_\infty, \Xx) = ((C_n/\Cc)^*, \cobdry)$.
The top-left submatrix $\ct{D}[1..i, 1..i]$ is the coboundary matrix for $\Cgr^*(X_\infty, X_{n-i})$. Indeed, this is precisely the full coboundary matrix with entries in $X_{n-i}$ removed.

It immediately follows that any procedure applied to matrix~$D$ that computes the intervals and generators of persistent absolute (resp.\ relative) homology, when applied to matrix~$\ct{D}$ will give us the intervals and generators of persistent relative (resp.\ absolute) cohomology.

\subsection{Persistence by matrix decomposition}
\label{sec:r-=-dv}

In~\cite{vineyards}, Cohen-Steiner et al.\ explain how to view the computation of persistent homology as a matrix decomposition problem, finding a factorization $D = RU$, where matrix $R$ is \textbf{reduced} and $U$ is invertible upper-triangular. Here we recap the relevant definitions.


For any matrix~$A$, we define $\low_A(j)$ to be the index of the lowest non-zero entry in the $j$-th column of $A$ (that is, the largest index $i$ such that $A[i,j] \ne 0$); it is undefined if the column is zero. We say that matrix $R$ is reduced if $\low_R$ is injective (over its domain of definition).

In what follows it is more convenient to look at the inverse of $U$, matrix $V = U^{-1}$. The decomposition $D = RU$ becomes $R = DV$. Whereas neither decomposition is unique, Cohen-Steiner et al.~\cite{vineyards} show that the map $\low_R$ is. It is precisely this map that gives the persistence pairing: a class born at the step~$g$ of the filtration dies at the step~$h$ iff $\low_R(h) = g$.

Suppose we have a decomposition $R = DV$. If the column $R[i] = 0$
(so $\low_R(i)$ is undefined), then the column $V[i]$ is a cycle, by definition. Furthermore, since $V$ is invertible upper-triangular, its diagonal entries are non-zero, and $V[i]$ is a cycle that does not exist in $X_{i-1}$, i.e.~it is exactly the cycle born at $\Hgr_*(X_i)$. Similarly, if  $R[j] \neq 0$, then it is the cycle that falls in the kernel of the map  $\Hgr_*(X_{j-1}) \to \Hgr_*(X_j)$, and $V[j]$ is the chain that appears in $X_j$ and has that cycle as its boundary.

\subsection{Homology generators}
\label{sec:homgen}

To relate the matrix discussion to the algebra of the previous section, we observe that the new basis of Theorem \ref{thm:persistence-main} appears in the matrices $R$
and $V$. The generator $\hat\sigma_f$ of the infinite interval $[a_f, \infty)$ is the column $V[f]$; the generator $\hat\sigma_g$ of the interval $[a_g,a_h)$ is the column $R[h]$; the chain $\hat\sigma_h$ that kills it is the column $V[h]$.

\begin{Example*}
\label{ex:abs-hom}
See Figure~\ref{fig:matrix-example} for the matrices $R,D,V$ of our running example. The map $\low_R$ gives the absolute homology persistence intervals $\Pers(\Hgr_*(\Ss)) = \{[1, \infty)_0, [2,3)_0, [4,5)_1, [6, \infty)_2\}$ where the subscript indicates
the dimension of the homology class. Columns of the matrices $V$ and $R$ give the cycles generating each of the intervals:
they are 
$\hat\ssx_1 = V[1] = \ssx_1$, 
$\hat\ssx_2 = R[3] = \ssx_1 - \ssx_2$, 
$\hat\ssx_4 = R[5] = \ssx_3 - \ssx_4$, and 
$\hat\ssx_6 = V[6] = \ssx_6 - \ssx_5$, respectively.

\begin{figure}
    \includegraphics{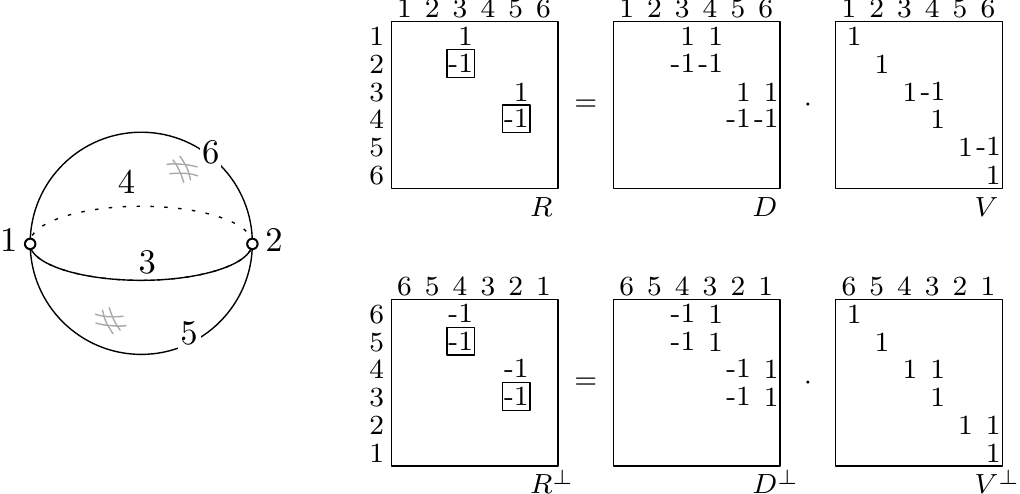}
    \caption{Filtration of a sphere, and the corresponding decompositions 
             $R = DV$ and $\ct{R} = \ct{D}\ct{V}$. Maps $\low_R$ and $\low_{\ct{R}}$ are
             shown with squares around the entries. We show only the non-zero 
             entries of the matrices.}
    \label{fig:matrix-example}            
\end{figure}
\end{Example*}

In view of Proposition \ref{prop:UCTkt}, we can also read off the generators for persistent relative homology: the column $V[f]$ is the generator for the interval $[-\infty, a_f)$, and the column $V[h]$ is the generator for $[a_g,a_h)$ (the chain that it represents becomes a relative cycle in $\Hgr_*(X_\infty, X_g)$, and remains nonzero until $\Hgr_*(X_\infty, X_h)$).

\begin{Example*}
\label{ex:rel-hom}
The four intervals of $\Pers(\Hgr_*(S_6, \Ss)) = \{[-\infty, 1)_0, [2,3)_1, [4,5)_2, [-\infty, 6)_2\}$ have respective generators 
$\hat\sigma_1 = V[1] = \ssx_1$, $\hat\sigma_3 = V[3] = \ssx_3$, 
$\hat\sigma_5 = V[5] = \ssx_5$, and $\hat\sigma_6 = V[6] = \ssx_6 - \ssx_5$.
\end{Example*}

\subsection{Cohomology generators}
\label{sec:cohogen}

We now take the global dual, and consider persistent relative and absolute cohomology. This time we compute the decomposition $\ct{R} = \ct{D} \ct{V}$ for the anti-transpose $\ct{D}$ of~ $D$.

We must take care to track the indices correctly. As a matrix, the rows and columns of $\ct{D}$ are labelled $\{1, \dots, n\}$ in the usual order. However, row~$i$ and column~$i$ refer to cell $\sigma_{n+1-i}$ in the original complex. If we define ${i}^* = n+1-i$, then we can think of the rows and columns as being labelled $\{n^*, \dots, 1^*\}$, so that row~$i^*$ and column~$i^*$ do indeed refer to cell $\sigma_i$. The numerical labels in Figure~\ref{fig:matrix-example} for $\ct{D}, \ct{R}, \ct{V}$ should be thought of as starred labels.

The columns of $\ct{R}$ and $\ct{V}$ contain the cocycles of $\Hgr^*(X_\infty, \Xx)$ and the cochains that kill them.
If $\low_{\ct{R}}(g^*) = h^*$, then there is a finite interval $[a_g,a_h)$ generated by the cocycle $\hat\sigma_h^* = \ct{R}[g^*]$ and killed by the cochain $\hat\sigma_g^* = \ct{V}[g^*]$.
If $\ct{R}[f^*] = 0$, then there is an infinite interval $[-\infty, a_f)$ generated by $\hat\sigma_f^* = \ct{V}[f^*]$.

\begin{Example*}
The persistent relative cohomology $\Pers(\Hgr^*(S_6, \Ss))$ has four intervals
$[-\infty,6)_2$, $[4,5)_2$, $[2,3)_1$, $[-\infty,1)_0$, 
generated respectively by
$\hat\ssx^*_6 = \ct{V}[6^*] = \ssx^*_6$, 
$\hat\ssx^*_5 = \ct{R}[4^*] = -\ssx^*_6 - \ssx^*_5$, 
$\hat\ssx^*_3 = \ct{R}[2^*] = -\ssx^*_3 - \ssx^*_2$, 
$\hat\ssx^*_1 = \ct{V}[1^*] = \ssx^*_1 + \ssx^*_2$.
\end{Example*}

Finally, we can read off the generators for persistent absolute cohomology $\Hgr^*(\Xx)$:
when $\ct{R}[f^*] = 0$,
the column $\ct{V}[f^*]$ is the cocycle which generates the interval $[-\infty, a_f)$; and when $\low_{\ct{R}(g^*)} = h^*$, the column $\ct{V}[g^*]$ is the cocycle which generates $[a_g,a_h)$.

\begin{Example*}
The persistent absolute cohomology $\Pers(\Hgr^*(\Ss))$ 
has four intervals 
$[6,+\infty)_2$, $[4,5)_1$, $[2,3)_0$, $[1,+\infty)_0$
generated respectively by
$\hat\ssx^*_6 = \ct{V}[6^*] = \ssx^*_6$, 
$\hat\ssx^*_4 = \ct{V}[4^*] = \ssx^*_4$, 
$\hat\ssx^*_2 = \ct{V}[2^*] = \ssx^*_2$,
$\hat\ssx^*_1 = \ct{V}[1^*] = \ssx^*_1 + \ssx^*_2$.
\end{Example*}

\subsection{Column and row algorithms}
\label{sec:elz}

The original persistence algorithm \cite{ELZ02, ZC05} finds the pairing
by processing matrix $D$ column-by-column to obtain the reduced matrix $R$.
In the context of $R = DV$ decomposition, one can express it as:
\begin{algorithm}
\caption{Column algorithm $\phcol$.}
\label{alg:column}
\begin{algorithmic}
    \STATE $R = D; V = I$
    \FOR{$i = 1$ to $n$}
        \WHILE{$\exists~j < i$ with $\low_R(j) = \low_R(i)$}
            \STATE $c = R[i][\low_R i]/R[j][\low_R j]$
            \STATE $R[i] = R[i] - cR[j]$
            \STATE $V[i] = V[i] - cR[j]$
        \ENDWHILE
    \ENDFOR
\end{algorithmic}
\end{algorithm}

\noindent
Here the definition of the constant~$c$
%
ensures that the
lowest non-zero element in column $i$ moves up after each iteration of the while
loop.
The condition of the while loop immediately implies that matrix $R$ is reduced
when the algorithm terminates. Furthermore, since we perform identical updates
on $R$ and $V$, we get an $R = DV$ decomposition.

The algorithm \phcol\ is essentially Gaussian elimination performed using column operations.
More commonly one would use column operations, processing the matrix row-by-row from the bottom up:
\begin{algorithm}
\caption{Row algorithm $\phrow$.}
\label{alg:row}
\begin{algorithmic}
    \STATE $R = D; V = I$
    \FOR{$i = n$ down to $1$}
        \STATE $\algvar{indices} = [ j \mid \low_R(j) = i ]$     \COMMENT{lows in the row $i$ of $R$}
        \STATE $p = \algvar{indices}[0]$                        \COMMENT{pivot}
        \FOR{$j \in \algvar{indices}[1..]$}
            \STATE $c = R[j][\low_R j]/R[p][\low_R p]$
            \STATE $R[j] = R[j] - cR[p]$
            \STATE $V[j] = V[j] - cV[p]$
        \ENDFOR
    \ENDFOR
\end{algorithmic}
\end{algorithm}

\noindent
It is not difficult to see that this algorithm also produces an $R = DV$
decomposition where matrix $R$ is reduced, and matrix $V$ is invertible
upper-triangular. 
What is less immediate is that the two algorithms produce identical
decompositions, so we prove this fact formally. (Notice that the statement would not be true if, during \phrow, we tried to cancel all non-zero elements in row~$i$ of~$R$, rather than restricting attention to the columns picked out by \textsf{indices}.) 

\begin{theorem}[Identical Output Theorem]
  The decompositions $R_c = DV_c$ and $R_r = DV_r$ produced by column and row
  algorithms respectively are identical, i.e.~$R_c = R_r$ and $V_c = V_r$.
\end{theorem}
\begin{proof}
We observe that once it determines the lowest non-zero element in a given column
of matrix $R$, neither algorithm changes that column in any subsequent
operations. Given a matrix $R = D$ we prove the claim by induction. The first column
with the lowest non-zero entry in $R$ is not modified by either algorithm.
Suppose that the columns with the lowest non-zero entries below $i$ are identical
in both $R_c$ and $R_r$, and $V_c$ and $V_r$. During the computation of the
column with the lowest non-zero entry in row $i$ we add columns with 
$\low_R > i$ in a decreasing order dictated by the lowest non-zero entry of the
column. Since the order and the columns are identical, so is the result.
\end{proof}

\begin{Remark}
Recently Milosavljevic, Morozov, and Skraba \cite{MMS11} showed that one can compute
persistence in matrix multiplication time.
\end{Remark}

\begin{Remark}
One can apply the two algorithms of this section to the restricted matrix $D_p$
that gives only the boundaries of the $p$-dimensional cells. We can still
extract some information from the $R_p = D_p V_p$ decomposition of this matrix:
the finite intervals $[g,h)$ in dimension $p-1$ and the births in dimension $p$,
i.e.~the endpoints $g$ or $h$ of any $p$-dimensional interval.
\end{Remark}



\section{Optimizations}
\label{sec:optimizations}

\subsection{Cohomology algorithm}
\label{sec:dmv}

One of our goals has been to relate our present work to an algorithm \pco\ for
persistent absolute cohomology that we described in~\cite{circular}. We based
that algorithm on the idea of maintaining a \emph{right filtration} (defined
in~\cite{zigzags}); as a result it looks different from $\phcol$ and $\phrow$ above.
In fact, we now show that one can view \pco\ as an optimization of \phrow\ applied to the matrix $\ct{D}$. We begin by reviewing the algorithm:
\begin{algorithm}
\caption{Cohomology algorithm $\pco$.}
\label{alg:cohomology}
\begin{algorithmic}
    \STATE {$\ct{Z} = [], \algvar{birth} = []$}
    \FOR{$i = 1$ to $n$}
        \STATE $\algvar{indices} = [ j ~|~ \ssx^*_i \in \cobdry z^*_j, z^*_j$ unmarked in $\ct{Z} ]$
        \IF{$\algvar{indices}$ are empty}
            \STATE prepend $\ssx^*_i$ to $\ct{Z}$ and $i$ to $\algvar{birth}$
        \ELSE
            \STATE prepend a marked $\ssx^*_i$ to $\ct{Z}$ and $i$ to $\algvar{birth}$
            \STATE $p = \algvar{indices}[0]$
            \FOR{$j = 1$ to $\operatorname{size}(\algvar{indices})$}
                \STATE $c = (\cobdry \ct{Z}[\algvar{indices}[j]])[i]/(\cobdry \ct{Z}[p])[i]$
                \STATE $\ct{Z}[\algvar{indices}[j]] = \ct{Z}[\algvar{indices}[j]] - c\ct{Z}[p]$
            \ENDFOR
            \STATE mark $\ct{Z}[p]$ and output the pair $[\algvar{birth}[p], i)$
        \ENDIF
    \ENDFOR
\end{algorithmic}
\end{algorithm}

\noindent
List $\ct{Z}$ maintains the cocycle basis for $H^*(X_i)$ in the right filtration order dictated by the
filtration of the space.  The marking above is for exposition only, in practice we drop a cocycle from the list $\ct{Z}$ as soon as it dies. When a new cell $\ssx_i$ enters, it is necessarily a cocycle (since it has no cofaces), but it may fall into a coboundary of a former cocycle, in which case (the else clause) we update the right filtration and drop the cocycle that $\ssx_i$ kills.

\begin{figure}
    \includegraphics{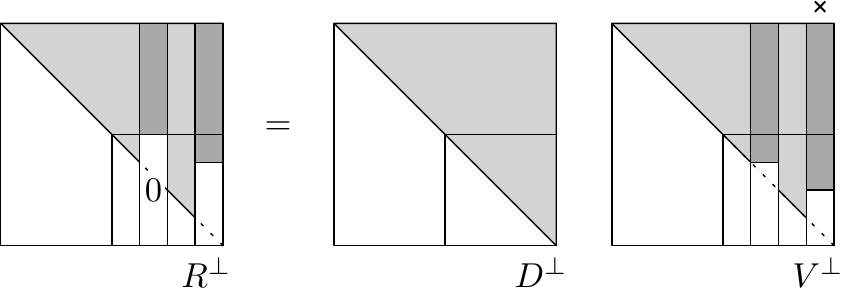}
    \caption{The structure of matrices $\ct{R} = \ct{D} \ct{V}$ during the execution of the row algorithm.}
    \label{fig:matrices}
\end{figure}

To see that this algorithm is a variation of the row algorithm from the previous section,
observe that the cocycles that it maintains are stored in the bottom-right corner of matrix $\ct{V}$
during the execution of the row algorithm.

\begin{claim}
The matrix $\ct{Z}$ in the cohomology algorithm after iteration $i$ is equal to the bottom-right corner
of the matrix $\ct{V}[(n-i)..n, (n-i)..n]$ after the $i$-th iteration of the row algorithm.
\end{claim}
\begin{proof}
We prove the claim inductively. Denoting with $\ct{R}_i, \ct{V}_i, \ct{Z}_i$ the various
matrices after $i$ iterations of both algorithms, assume the unmarked cocycles
$z^*_j$ in $\ct{Z}_i$ are exactly the cocycles with $\low_{\ct{R}}(j) > i$. In other
words, the corresponding columns $\ct{R}_i[j] = 0$. Furthermore assume that the two
matrices are identical, i.e.~$\ct{V}_i = \ct{Z}_i$. The claim is true when $i = 0$.
Our goal is to show it is true for $i = k$ assuming it is true for $i = k-1$.

At the $k$-th iteration, if cell $\ssx^*_k$ does not appear in the coboundary of
any cocycle, then its row in $\ct{R}_{k-1} = \cobdry \ct{V}_{k-1} = \cobdry \ct{Z}_{k-1}$ is zero. It
follows that it is not in the image of the map $\low_{\ct{R}_{k-1}}$ and therefore neither
algorithm performs any changes, so $\ct{V}_k = \ct{Z}_k$, and unmarked cocycles remain
as claimed.

If cell $\ssx^*_k$ is in the coboundary of a cocycle $z^*_j$ then $k$ is in the
image $\im \low_{\ct{R}_{k-1}}$. Moreover, from the inductive hypothesis the
indices $j$ of the columns of $\ct{R}_{k-1}$ that have $\low_{\ct{R}_{k-1}}(j) = i$
are exactly the unmarked cocycles in $\ct{Z}_{k-1}$ that have $\ssx^*_k$ in their
coboundary. Therefore, the update performed by both algorithms is identical.
\end{proof}

\begin{Remark}
Since the matrix $R$ contains the final persistence pairing, expressed as the
map $\low_R$, the algorithm \phcol\ is commonly optimized to keep track only of
this matrix (and ignore matrix $V$). In contrast, \pco\ maintains only matrix $\ct{Z} = \ct{V}$.
\end{Remark}

\subsection{Practice}
\label{sec:practice}

The algorithm \pco\ above highlights the difference between the column and the
row versions of the persistence algorithm. \phcol\ stores all the dead cycles
since it has no choice: any of them might be required at some future point in the reduction. \phrow, on the
other hand, is able to `examine the future' by inspecting any chosen row. It is therefore free to drop a column once it has determined its pairing and used it in the update. \pco\ does so explicitly.

In practice, such row access may be difficult when computing homology: it requires quick access to the coboundary of a given cell (since that is what a row of $D$ is). In simplicial complex implementations it is common to represent simplices as lists of vertices; then their boundary maps are easy to compute on the fly, while their coboundaries require a full preprocessing of the entire boundary matrix. By switching to cohomology we turn the tables: all the primitives necessary for the row algorithm (and in particular the optimized version given in this section) are readily available.

\subsection{Experiments}
\label{sec:experiments}

The practical improvement resulting from these observations is startling. In the following table we compare the traditional persistent homology algorithm \phcol\ with the cohomology
algorithm \pco. We list the total number of operations performed (in terms of primitive operations during chain arithmetic), total running time, and peak space usage in terms of the number of elements stored.

\begin{center}
    \begin{tabular}{r|c|r|r|r}
        Dataset                & Algorithm   &       Operations        &   Time    &   Peak elements       \\
        \hline\hline
        \ds{M-50}              &  \pco\      &       2,171,909,275     &   106 s   &    575,758            \\
                               &  \phcol\    &     609,477,028,616     &  4160 s   &  6,461,866            \\
        \hline\hline 
        \ds{T-10,000}          & \pco\       &          55,930,317     &     6 s   &     22,629            \\
                               &  \phcol\    &      29,760,159,689     &   207 s   &    693,031
    \end{tabular}
\end{center}

We used the \Cpp library Dionysus \cite{dionysus} to perform the above experiments.
The homology algorithm $\phcol$ in the above table computes only the matrix~$R$
since it suffices to extract the barcode. It also uses the original optimization
of~\cite{ELZ02} and stores the non-zero coefficients only in the rows that correspond to the positive cells.
\ds{M-50} is a filtration of an 8-skeleton of a Rips complex built on 50 random points of 
a Mumford dataset~\cite{mumfordTDA,mumford_original} up to the maximum pairwise distance of 1.5; the largest complex consists of 663,901 simplices.  \ds{T-10,000} is an alpha shape filtration of 10,000 points sampled on a torus embedded in $\mathbf{R}^3$;  the size of the Delaunay triangulation is 557,727 simplices.
The speed-up is encouraging. We would like to point out that these examples are not cherry-picked: we have yet to find a filtration on which \phcol\ is the faster of the two.


\bigskip
{\bf Conclusion.} When combined, the algebraic and experimental observations
suggest that if given a choice, one is better off using the cohomology
algorithm. Most of the time one has such a choice: for example, when computing
only the persistence diagram.

\bibliographystyle{unsrt}
\bibliography{4c}









\end{document}